\documentclass{article}
\usepackage{texdraw}

\newtheorem{Thm}{Theorem}[section]
\newtheorem{Cor}{Corollary}[section]
\newtheorem{Prop}{Proposition}[section]

\newtheorem{Def}{Definition}[section]

\newcommand{\RR}{{\bf R}}
\newcommand{\ZZ}{{\bf Z}}

\newcommand{\qed}{\hspace*{\fill}Q.E.D.}  
\newcommand{\ProofEnd} {\hspace*{\fill}Q.E.D.}  

\title{\bf The Codimension One
  Homology
of \\
a Complete Manifold with\\
 Nonnegative Ricci Curvature}
\author{Zhongmin Shen and Christina Sormani}
\date{November 10, 2000}

\begin{document}
\maketitle

\begin{center}
{Abstract}
\end{center}
In this paper we prove that a complete noncompact
manifold with nonnegative Ricci curvature has a trivial
codimension one homology
unless it is a  flat normal bundle over a compact
totally geodesic submanifold.
In particular, we prove the conjecture that
a complete noncompact manifold with positive Ricci curvature
has a trivial codimension one integer homology.
We also
have a corollary stating when the codimension two integer homology of
such a manifold is torsion free.

\section{Introduction}

 In this paper we prove that a complete noncompact
manifold with nonnegative Ricci curvature has a trivial
codimension one homology
unless it is a split or   flat normal bundle over a compact
totally geodesic submanifold (Thms~\ref{ThmA} and~\ref{ThmB}).
In particular, we prove the conjecture that
a complete noncompact manifold with positive Ricci curvature
has a trivial codimension one integer homology (Cor~\ref{Cor1}).
We also
have a corollary stating when the codimension two integer homology of
such a manifold is torsion free (Cor~\ref{Cor2}).

In 1976, Shing Tung Yau proved that a complete noncompact
manifold with positive    Ricci curvature has a trivial
codimension one real    homology [Yau].  It was then conjectured that
the integer homology might be trivial as well.  In 1993, the
first author proved that if the manifold is also proper, then
the codimension one integer homology is trivial [Sh1].
Recall that a manifold is proper if the manifold's Busemann function   
has compact level sets.  It is unknown whether
all complete noncompact manifolds with positive Ricci curvature
are proper.  However, the authors have proven that such a manifold
is proper when it
has Euclidean volume growth [Sh2] or linear volume growth [So1].

Recently Itokawa and Kobayashi have proven that if a complete noncompact
manifold with nonnegative Ricci curvature has more than linear
volume growth then the codimension one integer homology is trivial
[ItKo, Thm 1].  Since a manifold with nonnegative Ricci curvature must 
have
at least linear volume growth [Yau],
the three papers [ItKo], [So1] and [Sh1] prove the conjecture
regarding manifolds with positive Ricci curvature.  In fact they     
demonstrate that codimension one integer homology is trivial unless
the manifold has linear volume growth and nonnegative Ricci curvature.
Note that split or flat normal bundles over a compact manifolds
with nonnegative Ricci curvature are examples of such
exceptional manifolds.
It is not hard to see that these examples can have
codimension one integer homology equal to $\ZZ$.

Itokawa and Kobayashi have also proven that if a manifold with 
nonnegative
Ricci curvature has bounded diameter growth then the integer
codimension one homology is either $\ZZ$, $\ZZ_2$ or $0$ [ItKo Thm 2].   
Note that the split or flat normal bundles over compact manifolds
have bounded diameter growth, but there also exist examples of manifolds
with nonnegative Ricci curvature, linear volume growth and
unbounded diameter growth [So1].
Itokawa and Kobayashi used minimizing currents to prove their results, 
and
the first author used an adaption of Morse theory to prove his result
regarding proper manifolds.  

We reprove these results and complete
the classification of the codimension one integer homology
using Poincare Duality, the Cheeger Gromoll Splitting Theorem
[ChGl] and properties of noncontractible loops developed by    
the second author [So2].  

\begin{Thm} \label{ThmA} \label{A}
Let $M^n$ be an orientable complete noncompact manifold with nonnegative
Ricci curvature and $G$ is an abelian group, then one of the following 
holds:
\begin{eqnarray*}
(i) & H_{n-1}(M, G) & = \,\,0,\\
(ii)& H_{n-1}(M, G) & =  \,\, G, \\
(iii)& H_{n-1}(M, G) & = \,\, \ker(G \stackrel {\times 2} 
{\rightarrow} G).
\end{eqnarray*}

Case (ii) can only occur when $M^n$ is an
an isometrically split manifold over a compact totally geodesic
orientable submanifold.

Case (iii) can only occur if $M^n$ is
a one-ended flat normal bundle over a compact totally geodesic
unorientable submanifold.

\end{Thm}

Recall that $\ker(\ZZ \stackrel {\times 2} {\rightarrow} \ZZ)=0$
and $\ker(\ZZ_k \stackrel {\times 2} {\rightarrow} \ZZ_k)=0$ if
$k$ is odd and is $\ZZ_2$ if $k$ is even.       

\begin{Thm} \label{ThmB} \label{B}
Let $M^n$ be a complete noncompact unorientable manifold with nonnegative Ricci
curvature and $G=\ZZ_2$ or $\ZZ$. Then one of
the following holds:
\begin{eqnarray*}
(i)  & H_{n-1}(M, G) & = \,\, 0, \\
(ii) & H_{n-1}(M,G)  & = \,\, G,  \\
(iii)& H_{n-1}(M, G) & = \,\,\ker(G \stackrel {\times 2} {\rightarrow} 
G).
\end{eqnarray*}     
Case (ii) can only occur when M is a one-ended flat normal bundle over a
compact totally geodesic orientable submanifold.

Case (iii) can only occur when M is an isometrically split manifold
over a compact unorientable submanifold.

Remark: These two cases are opposite to the two cases in 
Theorem~\ref{ThmA}.
\end{Thm}

For an outline of our proof of Theorems~\ref{A} and~\ref{ThmB}
and the various cases involved, see the following diagram:

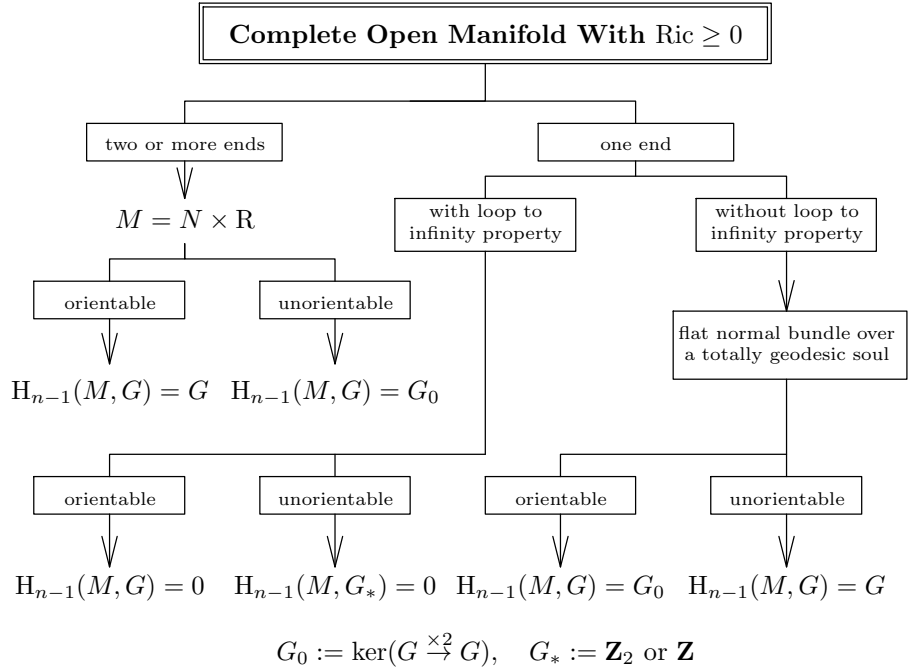
\begin{center}
\begin{texdraw}
\drawdim cm
\linewd  0.02

\move(1 0.3)
\textref h:C v:T \htext{\bf Complete Open Manifold With ${\rm Ric} \geq 
0$}

\move(-2.8 0.6) \lvec(4.8 0.6)\lvec(4.8 -0.2)\lvec(-2.8 -0.2)\lvec(-2.8 
0.6)
\move(-2.75 0.55) \lvec(4.75 0.55)\lvec(4.75 -0.15)\lvec(-2.75 
-0.15)\lvec(-2.75 0.55)

\move(1 -0.2)
\lvec(1 -0.7)
\lvec(-3 -0.7)
\lvec(-3 -1)
\move(-3 -1.2)
\textref h:C v:T \htext{\scriptsize{two or more ends}}

\move(-4.3 -1)
\lvec(-1.7 -1)
\lvec(-1.7 -1.5)
\lvec(-4.3 -1.5)
\lvec(-4.3 -1)

\move(0 -0.7)
\lvec(3 -0.7)
\lvec(3 -1)
\move(3 -1.2)
\textref h:C v:T \htext{\scriptsize{one end}}

\move(1.7 -1)\lvec(4.3 -1)\lvec(4.3 -1.5)\lvec(1.7 -1.5)\lvec(1.7 -1)

\move(-3 -1.5)
\arrowheadtype t:V \avec(-3 -2)

\move(-3 -2.3)
\textref h:C v:C \htext{$M=N\times {\rm R}$}

\move(-3 -2.6)
\lvec(-3 -2.8)
\lvec(-4 -2.8)
\lvec(-4 -3.1)
\move(-4 -3.4)
\textref h:C v:C \htext{\scriptsize{orientable}}
\move(-4 -3.6)
\arrowheadtype t:V \avec(-4 -4.2)
\move(-4 -4.6)
\textref h:C v:C \htext{${\rm H}_{n-1}(M, G) =G $}

\move(-5 -3.1)\lvec(-3 -3.1)\lvec(-3 -3.6)\lvec(-5 -3.6)\lvec(-5 -3.1)

\move(-3 -2.6)
\lvec(-3 -2.8)
\lvec(-1 -2.8)
\lvec(-1 -3.1)
\move(-1 -3.4)
\textref h:C v:C \htext{\scriptsize{unorientable}}

\move(-2 -3.1)\lvec(0 -3.1)\lvec(0 -3.6)\lvec(-2 -3.6)\lvec(-2 -3.1)

\move(-1 -3.6)
\arrowheadtype t:V \avec(-1 -4.2)
\move(-1 -4.6)
\textref h:C v:C \htext{${\rm H}_{n-1}(M, G)=G_0$}

\move(3 -1.5)
\lvec(3 -1.7)
\lvec(1 -1.7)
\lvec(1 -2)
\move(1 -2.2)
\textref h:C v:C \htext{\scriptsize{with loop to}}
\move(1 -2.5)
\textref h:C v:C \htext{\scriptsize{infinity property}}
\move(1 -2.7)
\lvec(1 -5.4)

\move(-0.2 -2)\lvec(2.2 -2)\lvec(2.2 -2.7)\lvec(-0.2 -2.7)\lvec(-0.2 -2)

\move(-1 -5.4)
\lvec(-1 -5.7)
\move(-1 -6)
\textref h:C v:C \htext{\scriptsize{unorientable}}
\move(-1 -6.2)
\arrowheadtype t:V \avec(-1 -6.9)
\move(-1 -7.2)
\textref h:C v:C \htext{${\rm H}_{n-1}(M,G_*)=0$}

\move(-2 -5.7)\lvec(0 -5.7)\lvec(0 -6.2)\lvec(-2 -6.2)\lvec(-2 -5.7)

\move(1 -5.4)
\lvec(-4 -5.4)
\lvec(-4 -5.7)
\move(-4 -6)
\textref h:C v:C \htext{\scriptsize{orientable}}
\move(-4 -6.2)
\arrowheadtype t:V \avec(-4 -6.9)
\move(-4 -7.2)
\textref h:C v:C \htext{${\rm H}_{n-1}(M, G) = 0$}

\move(-5 -5.7)\lvec(-3 -5.7)\lvec(-3 -6.2)\lvec(-5 -6.2)\lvec(-5 -5.7)

\move(3 -1.7)
\lvec(5 -1.7)
\lvec(5 -2)
\move(5 -2.2)
\textref h:C v:C \htext{\scriptsize{without loop to}}
\move(5 -2.5)
\textref h:C v:C \htext{\scriptsize{infinity property}}
\move(3.8 -2)\lvec(6.2 -2)\lvec(6.2 -2.7)\lvec(3.8 -2.7)\lvec(3.8 -2)

\move(5 -2.7)
\arrowheadtype t:V \avec(5 -3.5)
\move(5 -3.8)
\textref h:C v:C \htext{\scriptsize{flat normal bundle over}}
\move(5 -4.1)
\textref h:C v:C \htext{\scriptsize{a totally geodesic soul}}

\move(3.5 -3.5)\lvec(6.6 -3.5)\lvec(6.6 -4.4)\lvec(3.5 -4.4)\lvec(3.5 
-3.5)

\move(5 -5.4)
\lvec(2 -5.4)
\lvec(2 -5.7)
\move(2 -6)
\textref h:C v:C \htext{\scriptsize{orientable}}
\move(2 -6.2)
\arrowheadtype t:V \avec(2 -6.9)
\move(2 -7.2)
\textref h:C v:C \htext{${\rm H}_{n-1}(M, G) =G_0$}

\move(1 -5.7)\lvec(3 -5.7)\lvec(3 -6.2)\lvec(1 -6.2)\lvec(1 -5.7)

\move(5 -4.4)
\lvec(5 -5.7)
\move(5 -6)
\textref h:C v:C \htext{\scriptsize{unorientable}}
\move(5 -6.2)
\arrowheadtype t:V \avec(5 -6.9)
\move(5 -7.2)
\textref h:C v:C \htext{${\rm H}_{n-1}(M, G) =G$}

\move(4 -5.7)\lvec(6 -5.7)\lvec(6 -6.2)\lvec(4 -6.2)\lvec(4 -5.7)

\move(1 -8)
\textref h:C v:C \htext{$ G_0 := \ker(G \stackrel {\times 2} 
{\rightarrow} G), \ \ \ G_* := {\bf Z}_2 \ {\rm or} \ {\bf Z}$}

\move(-6 0)
\lcir r:0.0
\move(8 0)
\lcir r:0

\move(0 1)
\move(0 -10)
\lcir r:0

\end{texdraw}
\end{center}

\begin{Cor} \label{Cor1}
Let $M^n$ be a complete noncompact manifold with nonnegative Ricci  
curvature.  If there is a point $p\in M$
such that $Ricci_p(v,v)>0$ then $H_{n-1}(M, \ZZ)$ is trivial.
\end{Cor}

Note that we have eliminated the possibility that $H_{n-1}(M, \ZZ)= 
\ZZ_2$.
This was indicated as a possibility in [ItKo, Thm 2] if $M$ had
a cover of degree one or two isometric to
$\bar{M}\times [0,\infty) \cup K$ where $K$ and $\bar{M}$ are compact.
Their second theorem also indicated that $\ZZ$ was a possibility
only if $M$ had a cover of degree one or two which split isometrically
exactly as we have.

The paper is divided into two sections.
In the first section,
we define a topological property called the {\em loops to infinity 
property}
[Defn~\ref{looptoinf}].  We then
prove that orientable manifolds with this property and only one end
have trivial codimension
one $G$ homology where $G$ is any abelian group
[Prop ~\ref{PartI}].  This section
of the paper is purely topological and uses Poincare Lefschetz Duality 
and
the Universal Coefficient Theorem (c.f. [Mun] and [Mas]).

In the second section of the paper, we eliminate the
topological conditions in Propositions~\ref{PartI}
using the properties of manifolds with nonnegative
Ricci curvature [Props ~\ref{PartII}-~\ref{PartIV}].
In particular, we use the Splitting Theorem [ChGl] and
a theorem in [So2].

In addition to proving Theorems~\ref{A} and~\ref{ThmB}, we prove the
following corollary in Section 2.

\begin{Cor} \label{Cor2}   
Let $M^n$ be an orientable complete noncompact manifold with nonnegative
Ricci curvature, then $M$ is either an isometrically split manifold
or a flat normal bundle over a compact totally geodesic submanifold
or, for any abelian group $G$,
$$
H_{n-2}(M, \ZZ) \ast G=0,
$$
in which case $H_{n-2}(M, \ZZ)$ has no elements of finite order.
\end{Cor}


The authors would like to thank Dan Christensen for helpful discussions
regarding homology theory.

\section{Topology and $H_{n-1}(M, \ZZ)$}

Before we prove Theorems~\ref{A} and~\ref{B}
we introduce some topological
concepts from
[So2] concerning noncontractible loops in complete noncompact manifolds.

\begin{Def} \label{LoopDefn} {\em
Given a ray $\gamma$ and a loop $C:[0,L] \to M$ based at $\gamma(0)$,
we say that a loop $\bar{C}:[0,L] \to M$ is
{\em homotopic to} $C$ along $\gamma$
if there exists $r>0$ with $\bar{C}(0)= \bar{C}(L)= \gamma(r)$ and 
the
loop, 
constructed by joining $\gamma$ from $0$ to $r$ with $\bar{C}$ from $0$
to $L$ and then with $\gamma$ from $r$ to $0$, is homotopic to $C$
(see diagram below). }
\end{Def}

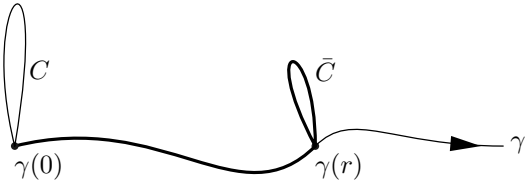
\begin{center}
\begin{texdraw}
\drawdim cm
\linewd  0.02


\move(-2 0)
\linewd 0.04
\clvec(0 0.5)(1 -1)(2 0)
\fcir f:0.1 r:0.05
\linewd 0.02
\clvec(2.5 0.5)(3 0)(4.5 0)
\move(4 0.02)
\arrowheadtype t:F \avec(4.2 0)

\move(-2 0)
\fcir f:0.1 r:0.05
\clvec(-1.5 3)(-2.5 2)(-2 0)


\move(2 0)
\linewd 0.04
\clvec(1.2 1.5)(2 1.5)(2 0)

\move(2 1)
\textref h:L v:C \htext{$\bar{C}$}

\move(-1.8 1)
\textref h:L v:C \htext{$C$}

\move(4.6 0)
\textref h:L v:C \htext{$\gamma$}
\move(2 -0.1)
\textref h:L v:T \htext{$\gamma(r)$}
\move(-2 -0.1)
\textref h:L v:T \htext{$\gamma(0)$}

\move(0 2.5)
\fcir f:1 r:0.0
\move(0 -1)
\fcir f:1 r:0.0
\move(-5.5 0)
\fcir f:1 r:0.1
\move(5.5 0)
\fcir f:1 r:0.1
\end{texdraw}
\end{center}

\begin{Def} \label{looptoinf} {\em
A complete noncompact manifold $M^n$ has the {\em loops to infinity
property} if there is a ray $\gamma$, such that for any
element, $g\in \pi_1(M, \gamma(0))$, and any compact set,
$K \subset M$, there exists a closed loop, $\bar{C}$,
contained in $M\backslash K$ which is
homotopic along $\gamma$ to a representative loop, $C$,
such that $g= [C]$.

When a ray, $\gamma$, is specified, we say $M^n$ has the loops
to infinity property along $\gamma$.}
\end{Def}

In [So2] 
it is proven that a manifold with nonnegative Ricci
curvature either has the loops to infinity property along some ray
or the manifold is a flat normal bundle over a compact totally
geodesic soul and has a double cover which splits isometrically.

Recall that a manifold is said to have {\em one end} if given any
compact set $K \subset M$, there is only one unbounded connected
component of $M \setminus K$.  Note that if a complete manifold
has two or more ends then it has a length minimizing geodesic
defined for all values of $\RR$, called a {\em line}.  Recall also
that by the Cheeger-Gromoll Splitting Theorem [ChGl], if a complete
noncompact manifold  with nonnegative Ricci curvature contains a line,
then it splits isometrically.  So a manifold with nonnegative Ricci
curvature either splits isometrically or has only one end.

Thus, in some sense, most manifolds with nonnegative Ricci curvature
have the loops to infinity property and only one end. 
We now prove the following simplified version of Theorem~\ref{ThmA}
relating the loops to infinity property to the codimension one integer
homology of a manifold.

\begin{Prop} \label{PartI}
Let $M^n$ be a complete noncompact manifold with one end and
the loops to infinity property.  Then
$H_{n-1}(M,\ZZ_2)$
is trivial. 

If $M^n$ is also orientable, then $H_{n-1}(M,G)$ is trivial
where $G$ is an abelian group.
\end{Prop}

\bigskip   \noindent {\bf Proof:}  
Suppose $H_{n-1}(M,G)$ is not trivial.  Then there is a chain, $\sigma$,
which has no boundary but which is not
a boundary .  See the diagram below.  The images of the simplices of 
this
chain must be contained in some compact set, $K_0= Cl(U_0)$,
containing $\gamma(0)$.  Since $M^n$ has
only one end, there is a compact set, $K$, consisting of $K_0$ and the
bounded connected components of $M \setminus K_0$.  Furthermore, $K$
can be chosen with smooth boundary.  Thus $K$ is an n dimensional
manifold with boundary such that $\partial K$ is connected and
\begin{equation} \label{contra}
H_{n-1}(K, G) \neq 0.
\end{equation}

\begin{center}
\begin{texdraw}
\drawdim cm
\linewd  0.02



\move(0 1)
\clvec(0.2 0.7)(0.2 -0.7)(0 -1)
\lfill f:0.95
\lpatt(0.06 0.08)
\clvec(-0.2 -0.7)(-0.2 0.7)(0 1)
\lfill f:0.95
\lpatt()


\move(-2 -0.8)
\fcir f:0.1 r:0.05
\clvec(-1 -0.3)(-0.5 -0.5)(0 -0.5)
\fcir f:0.1 r:0.05
\clvec(0 -0.5)(1 -0.7)(3 -0.5)
\move(2.8 -0.52)
\arrowheadtype t:F \avec(3 -0.5)


\move(0 1)
\clvec(1 1.2)(1.5 1.2)(2.5 1.5)
\move(0 1)
\clvec(-1 0.8)(-3 1)(-3 -0.7)
\clvec(-3 -1.5)(-1 -0.8)(0 -1)
\clvec(1 -1.2)(2.5 -1.3)(3 -1.5)


\move(-2.3 -0.2)
\lvec(-1.8 -0.5)\lvec(-1.5 -0.2)\lvec(-2.3 -0.2)
\lfill f:0.95

\move(-1.8 -0.5)
\lvec(-1 0)\lvec(-1.5 -0.2)\lvec(-1.8 -0.5)
\lfill f:0.8

\move(-1.5 -0.2)
\lvec(-1.2 0.4)\lvec(-1 0)\lvec(-1.5 -0.2)
\lfill f:0.92

\move(-1.5 -0.2)
\lvec(-2.3 -0.2)\lvec(-1.2 0.4)\lvec(-1.5 -0.2)
\lfill f:0.94

\move(-2.3 -0.2)
\lvec(-1.2 0.4)\lvec(-1.7 0.3)\lvec(-2.3 -0.2)
\lfill f:0.92

\move(-1.7 0.3)
\lpatt(0.04 0.06)
\lvec(-1.8 -0.5)
\move(-1.7 0.3)
\lvec(-1 0)


\move(0.2 -0.65)
\textref h:L v:T \htext{$\gamma(t_0)$}
\move(-2.1 -0.8)
\textref h:R v:C \htext{$\gamma(0)$}

\move(-2 0.8)
\textref h:C v:B \htext{$K$}
\move(0 1.1)
\textref h:C v:B \htext{$\partial K$}

\move(-1 0.2)
\textref h:L v:C \htext{$\sigma$}

\move(0 2.3)
\fcir f:1 r:0.0
\move(0 -2)
\fcir f:1 r:0.0
\move(-5.5 0)
\fcir f:1 r:0.1
\move(5.5 0)
\fcir f:1 r:0.1

\end{texdraw}
\end{center}

We now proceed to show that (\ref{contra}) cannot hold under
the appropriate conditions on $M$.

Let $\gamma$ be the ray along which $M^n$ has the loops
to infinity property.
Since $\gamma$ is a ray, it must leave $K$, so let
$t_0>0$ be defined such that $\gamma(t_0)\in \partial K$.
In [So2], it is proven that if $M^n$ has the loops to
infinity property along $\gamma$ and $\partial K$ is smooth and
connected,
then
\begin{equation}
i_*:\pi_1(\partial K, \gamma(t_0)) \to \pi_1(K, \gamma(t_0))
\end{equation}
is onto.  Here $i_*$ is the map induced by the inclusion map
$i: \partial K \to K$.

Now $H_1(N, \ZZ)$ is just the abelianization of $\pi_1(N,p)$
(c.f. [Mas, Ch. VIII, Thm 7.1]).  Thus we
know
\begin{equation}
i_*:H_1(\partial K, \ZZ) \to H_1(K, \ZZ)
\end{equation}
is onto.
This is part of a long exact sequence involving the relative homology
and the reduced homology,
\begin{equation}
H_1(\partial K, \ZZ) \stackrel {i_*}{\rightarrow} H_1(K, \ZZ)
 \stackrel{\pi_*}{\rightarrow} H_1(K, \partial K, \ZZ)
 \stackrel{\partial_*}{\rightarrow} \tilde{H}_0(\partial K, \ZZ),
\end{equation}
(c.f. [Mun, Thm 23.3]).  So $Ker(\pi_*)= Im(i_*)=H_1(K, \ZZ)$.  Thus
$\pi_*= 0$ and $\partial_*$ is one-to-one.  Now, $\partial K$ is
connected so the reduced homology, $\tilde{H}_0(\partial K, \ZZ)$, is
trivial.  Thus the relative homology, $H_1(K, \partial K, \ZZ)$, is 
trivial
as well.

We now claim that the cohomology $H^1(K, \partial K, G)$ must be
trivial as well.  To prove this we will apply the
Universal Coefficient Theorem for Cohomology
(c.f. [Mun, Thm 53.1]), which provides the following
short exact sequence:
$$   
0 \to Ext(H_0(K, \partial K; \ZZ), G) \to
H^1(K, \partial K; G) \to Hom(H_1(K, \partial K; \ZZ), G)
\to 0.
$$
Since $H_1(K, \partial K; \ZZ)$ is trivial, the homomorphisms
from it to $G$ must be trivial as well.  Thus
\begin{equation} \label{ExtH1}
Ext(H_0(K, \partial K; \ZZ), G) = H^1(K, \partial K; G).
\end{equation}
On the other hand, since
$K$ and $\partial K$ are connected, the long exact relative homology
sequence (c.f. [Mun, Thm 23.3]),
\begin{equation}
... \to H_0(\partial K,\ZZ) \stackrel{i_*}{\rightarrow}
H_0(K, \ZZ) \stackrel {\pi_*} {\rightarrow}
H_0(K, \partial K; \ZZ) \to 0
\end{equation}
is
\begin{equation}
... \to \ZZ \stackrel{i_*}{\rightarrow}
 \ZZ \stackrel{\pi_*} {\rightarrow}
H_0(K, \partial K; \ZZ) \to 0,
\end{equation}
where $i_*$ is an isomorphism.  Thus $ker(\pi_*)= Im(i_*)= \ZZ$
and $\pi_*$ is onto, so $H_0(K, \partial K; \ZZ)= 0$.   Thus
$Ext(H_0(K, \partial K; \ZZ), G)=0$ and, by (\ref{ExtH1}),
$H^1(K, \partial K, G)$ must be trivial as well.

Now by Poincare Lefschetz Duality [Mun, Thm 70.7], we know that
\begin{equation}
H^1(K, \partial K, \ZZ_2) \to H_{n-1}(K, \ZZ_2)
\end{equation}
is an isomorphism.   Since $H^1(K, \partial K, \ZZ_2)$ is trivial
and, by (\ref{contra}), $H_{n-1}(K, \ZZ_2)$ is not trivial, we have
a contradiction.  Thus $H_{n-1}(M, \ZZ_2)$ must be trivial.

If $(K, \partial K)$ is orientable, then Poincare Lefschetz Duality
implies that
\begin{equation}
H^1(K, \partial K, G) \to H_{n-1}(K, G)
\end{equation}
is an isomorphism.  Thus, when $M$ is orientable,
we contradict (\ref{contra}), so
$H_{n-1}(M, G)$ must be trivial.
\qed
\bigskip

\section{Ricci Curvature and $H_{n-1}(M, \ZZ)$}

To prove the main theorems, we must deal with complete noncompact
manifolds with nonnegative Ricci curvature that either
have more than one end, are not orientable or fail to satisfy the
loops to infinity property.  We deal with these three cases seperately. 
Refer to the diagram at the end of the paper.

\begin{Prop} \label{PartII}
Let $M^n$ be a complete noncompact with nonnegative Ricci
curvature and two or more ends.  Let $G$ be an abelian group.
Then $M= N^{n-1}\times \RR$ and
$$
H_{n-1}(M,G)= \cases{G,&if  M is orientable \cr
 \ker(G \stackrel {\times 2} {\rightarrow} G) ,& if M is not orientable 
}.
$$
\end{Prop}

\begin{Prop} \label{PartIII}
Let $M^n$ be a one-ended complete noncompact manifold
with nonnegative Ricci curvature satisfying the loops to infinity
property.  Then, regardless of the orientability of $M$,
$H_{n-1}(M,\ZZ)=0$.
\end{Prop}

\begin{Prop} \label{PartIV}
Let $M^n$ be a complete noncompact manifold with  nonnegative Ricci
curvature, that has one end and
doesn't have a ray with the loops to infinity
property. Let $G$ be an abelian group.
Then $M$ is a flat normal bundle over a totally
geodesic soul and
$$
H_{n-1}(M,G)= \cases{G,&if  M is unorientable \cr
 \ker(G \stackrel {\times 2} {\rightarrow} G) ,& if M is orientable }.
$$
\end{Prop}

Note that Propositions~\ref{PartI}, ~\ref{PartII}
and ~\ref{PartIV} imply Theorem~\ref{A}, and
Propositions~\ref{PartI}-~\ref{PartIV} imply Theorem~\ref{B}.
The proof of Corollary~\ref{Cor1} will appear at the end
of the paper.

\bigskip  
\noindent {\bf Proof of Proposition~\ref{PartII}:}
By the Cheeger-Gromoll Splitting Theorem, a complete noncompact
manifold with nonnegative Ricci curvature and two or more ends
must split isometrically, $M^n= N^{n-1}\times \RR$, where $N^{n-1}$
is compact.  Thus $H_{n-1}(M^n, G)= H_{n-1}(N^{n-1},G)$, which
is
$$
\ker(G \stackrel {\times 2} {\rightarrow} G)
$$
if $N^{n-1}$ is not orientable and is $G$ if $N^{n-1}$ is
orientable (c.f. [Mun Cor 65.5]).  The proposition then
follows from the fact that $N^{n-1}$ is  orientable iff $M^n$ is orientable.
\qed

\bigskip  
\noindent {\bf Proof of Proposition~\ref{PartIII}:}
This proposition holds if $M$ is orientable by Proposition~\ref{PartI}.
If $M$ is not orientable, then it has a double cover
$\pi:\tilde{M}\to M$ such that $\tilde{M}$
is orientable.   We first claim that $\tilde{M}$ has only one end.

Assume on the contrary that $\tilde{M}$ has two or more ends.
By the Cheeger-Gromoll Splitting Theorem,  it splits isometrically,
$\tilde{M}= N^{n-1}\times \RR$ where $N$ is compact [ChGl]. 
Since $\tilde{M}$ is orientable,
so is the totally geodesic submanifold, $N^{n-1}$.  Let $g$ be
the deck transform acting on $\tilde{M}$.  Then $g$ takes
lines to lines and acts as an isometry on each component,
$g:\RR \to \RR$ and $g:N^{n-1} \to N^{n-1}$.  Now
$g:\RR \to \RR$ is an isometry so $g(r)= \pm r +r_0$. 
If $g(r)= r +r_0$, then
$
\pi: \{p_0\}\times \RR \to \pi(\{p_0\}\times \RR)
$
would map a line to a line and $M$ itself would split
isometrically and have two ends.   Thus $g(r)= -r +r_0$.

Now $M$ has the loops to infinity property along a ray, $\gamma$,
which lifts to a ray $\tilde{\gamma}$. Since $N$ is compact, we can 
choose
$\tilde{\gamma}(t)=(q_0,  t+t_0)$.  Let $C_0$ be a
curve in $\tilde{M}$ joining $(q_0, t_0)$ to $g(q_0,t_0)=(gq_0, 
-t_0+r_0)$.  Then
for any compact set $K \subset M$, including $K=\pi(N^{n-1}\times 
\{r_0/2\})$,
there exists a curve, $C_1 \in M \setminus K$,
based at $\gamma(t_1)$ which is homotopic along $\gamma$ to
$\pi \circ C_0$ [Defn \ref{looptoinf}].  See the diagram below.
  By lifting the homotopy, we lift
$C_1$ to a curve, $\tilde{C}_1$, which runs from
$\tilde{\gamma}(t_1)=(q_0, t_1+t_0)$
to $g \tilde{\gamma}(t_1)=(g q_0, -t_1 -t_0 + r_0)$ [Defn 
\ref{LoopDefn}].
Thus $\tilde{C}_1$ passes through $N \times \{r_0/2\}$, and we have
a contradiction.  Any unorientable manifold with nonnegative Ricci
curvature and the loops to infinity property, has an orientable double
cover with only one end.

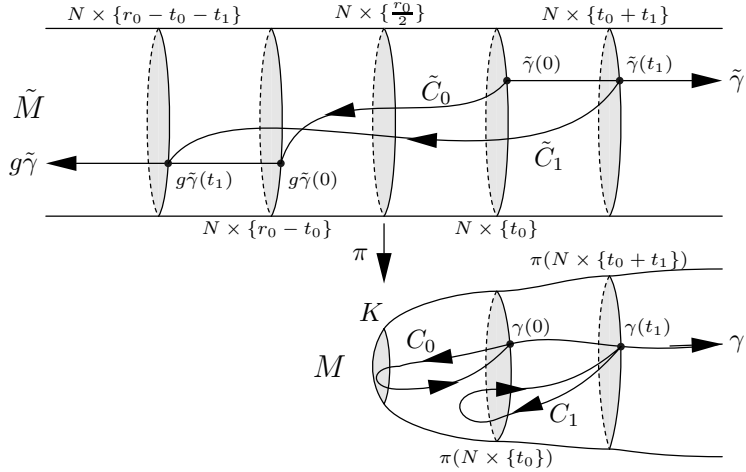
\vspace{.5cm}

\begin{center}
\begin{texdraw}
\drawdim cm
\linewd  0.02



\move(-4.5 3)
\lvec(4.5 3)
\move(-4.5 0.5)
\lvec(4.5 0.5)


\move(0 3)
\lpatt(0.05 0.06)
\clvec(-0.2 2.7)(-0.2 0.8)(0 0.5)
\lfill f:0.9
\lpatt()
\clvec(0.2 0.8)(0.2 2.7)(0 3)
\lfill f:0.9

\move(1.5 3)
\lpatt(0.05 0.06)
\clvec(1.3 2.7)(1.3 0.8)(1.5 0.5)
\lfill f:0.9
\lpatt()
\clvec(1.7 0.8)(1.7 2.7)(1.5 3)
\lfill f:0.9

\move(3 3)
\lpatt(0.05 0.06)
\clvec(2.8 2.7)(2.8 0.8)(3 0.5)
\lfill f:0.9
\lpatt()
\clvec(3.2 0.8)(3.2 2.7)(3 3)
\lfill f:0.9

\move(-1.5 3)
\lpatt(0.05 0.06)
\clvec(-1.7 2.7)(-1.7 0.8)(-1.5 0.5)
\lfill f:0.9
\lpatt()
\clvec(-1.3 0.8)(-1.3 2.7)(-1.5 3)
\lfill f:0.9

\move(-3 3)
\lpatt(0.05 0.06)
\clvec(-3.2 2.7)(-3.2 0.8)(-3 0.5)
\lfill f:0.9
\lpatt()
\clvec(-2.8 0.8)(-2.8 2.7)(-3 3)
\lfill f:0.9

\move(1.63 2.3)
\fcir f:0.1 r:0.05
\arrowheadtype t:F \avec(4.5 2.3)
\move(3.13 2.3)
\fcir f:0.1 r:0.05

\move(4.6 2.3)
\textref h:L v:C \htext{$\tilde{\gamma}$}

\move(-1.37 1.2)
\fcir f:0.1 r:0.05
\arrowheadtype t:F \avec(-4.5 1.2)
\move(-2.87 1.2)
\fcir f:0.1 r:0.05
\move(-4.6 1.2)
\textref h:R v:C \htext{$g\tilde{\gamma}$}

\move(3.13 2.3)
\clvec(2 0.5)(-2 2.5)(-2.87 1.2)
\move(0.5 1.515)
\arrowheadtype t:F \avec(0.3 1.52)

\move(1.63 2.3)
\clvec(1 1.5)(-1 2.5)(-1.37 1.2)
\move(-0.5 1.9)
\arrowheadtype t:F \avec(-0.8 1.848)


\move(0 0.4)
\arrowheadtype t:F \avec(0 -0.4)
\move(-0.2 0)
\textref h:R v:C \htext{$\pi$}



\move(0 -1)
\clvec(-0.2 -1.3)(-0.2 -1.7)(0 -2)
\clvec(0.1 -1.8)(0.1 -1.2)(0 -1)
\lfill f:0.9
\move(1.5 -0.5)
\lpatt(0.05 0.06)
\clvec(1.3 -0.8)(1.3 -2.2)(1.5 -2.5)
\lfill f:0.9
\lpatt()
\clvec(1.8 -2.2)(1.7 -0.7)(1.5 -0.5)
\lfill f:0.9
\move(3 -0.3)
\lpatt(0.05 0.06)
\clvec(2.8 -0.6)(2.8 -2)(3 -2.6)
\lfill f:0.9
\lpatt()
\clvec(3.2 -2)(3.2 -0.6)(3 -0.3)
\lfill f:0.9


\move(0 -1)
\clvec(0.2 -0.7)(1 -0.5)(1.5 -0.5)
\clvec(2 -0.5)(2.5 -0.3)(3 -0.3)
\clvec(3.2 -0.3)(3.5 -0.2)(4.5 -0.2)

\move(0 -2)
\clvec(0.2 -2.3)(1 -2.5)(1.5 -2.5)
\clvec(2 -2.5)(2.5 -2.6)(3 -2.6)
\clvec(3.2 -2.6)(3.5 -2.7)(4.5 -2.7)


\move(1.68 -1.2)
\fcir f:0.1 r:0.05
\clvec(2.5 -1)(3 -1.4)(4.5 -1.2)
\move(3.8 -1.225)
\arrowheadtype t:F \avec(4.5 -1.2)
\move(3.15 -1.24)
\fcir f:0.1 r:0.05


\move(3.15 -1.24)
\clvec(2.5 -1.8)(2 -1.8)(1.5 -1.8)
\clvec(0.7 -1.8)(1 -2.4)(1.5 -2.2)
\clvec(2 -2)(2.5 -2)(3.15 -1.24)

\move(1.7 -1.8)
\arrowheadtype t:F \avec(1.9 -1.8)
\move(2.2 -2)
\textref h:L v:T \htext{$C_1$}
\move(1.9 -2.05)
\arrowheadtype t:F \avec(1.7 -2.12)


\move(1.68 -1.2)
\clvec(0.5 -1.5)(0.5 -1.4)(0.2 -1.5)
\clvec(-0.2 -1.5)(-0.2 -1.8)(0.2 -1.8)
\clvec(0.8 -1.8)(1.2 -1.7)(1.68 -1.2)

\move(0.7 -1.4)
\arrowheadtype t:F \avec(0.5 -1.43)
\move(0.5 -1.3)
\textref h:C v:B \htext{$C_0$}
\move(0.8 -1.74)
\arrowheadtype t:F \avec(1 -1.708)


\move(4.6 -1.2)
\textref h:L v:C \htext{$\gamma$}
\move(3.2 -1.1)
\textref h:L v:B \htext{{\scriptsize $\gamma(t_1)$}}
\move(1.7 -1.1)
\textref h:L v:B \htext{{\scriptsize $\gamma(0)$}}

\move(1.5 -2.6)
\textref h:C v:T \htext{{\scriptsize $\pi(N \times \{t_0\})$}}

\move(3 -0.2)
\textref h:C v:B \htext{{\scriptsize $\pi(N \times \{ t_0+t_1\} )$}}

\move(0 -0.9)
\textref h:R v:B \htext{$K$}

\move(2.2 1.5)
\textref h:C v:T \htext{$\tilde{C}_1$}
\move(0.7 2)
\textref h:C v:B \htext{$\tilde{C}_0$}

\move(-2.77 1.1)
\textref h:L v:T \htext{{\scriptsize $g\tilde{\gamma}(t_1)$}}
\move(-1.27 1.1)
\textref h:L v:T \htext{{\scriptsize $g\tilde{\gamma}(0)$}}

\move(1.73 2.4)
\textref h:L v:B \htext{ {\scriptsize $\tilde{\gamma}(0)$}}
\move(3.23 2.4)
\textref h:L v:B \htext{{\scriptsize $\tilde{\gamma}(t_1)$}}

\move(3 3.05)
\textref h:C v:B \htext{{\scriptsize $N \times\{t_0+t_1\} $ }}
\move(1.5 0.45)
\textref h:C v:T \htext{{\scriptsize $N\times \{t_0\}$}}
\move(0 3.05)
\textref h:C v:B \htext{{\scriptsize $N \times\{{r_0\over 2}\} $ }}
\move(-1.5 0.45)
\textref h:C v:T \htext{{\scriptsize $N \times\{r_0-t_0\} $ }}
\move(-3 3.05)
\textref h:C v:B \htext{{\scriptsize $N \times\{r_0-t_0-t_1\} $ }}
\move(-4.5 2)
\textref h:R v:C \htext{{\large $\tilde{M}$}}
\move(-0.5 -1.5)
\textref h:R v:C \htext{{\large $M$}}

\move(0 3)
\fcir f:1 r:0.0
\move(0 -3)
\fcir f:1 r:0.0
\move(-5.5 0)
\fcir f:1 r:0.1
\move(5.5 0)
\fcir f:1 r:0.1

\end{texdraw}
\end{center}
  
Now $\tilde{M}$ has only one end and, if we prove that
$\tilde{M}$ has the loops to infinity property, then  it
satisfies all the conditions of Proposition~\ref{PartI}.

The ray, $\gamma$, in $M$ with the loops to
infinity property can be lifted to a ray $\tilde {\gamma}$ in
$\tilde{M}$. 
We claim that $\tilde{M}$ satisfies the loops to infinity property
along $\tilde{\gamma}$ [Defn~\ref{looptoinf}].  Given any
element, $[C] \in \pi_1(\tilde{M}, \tilde{\gamma}(0))$, $\pi \circ C$ is 

a loop in $M$ based at $\gamma(0)$.  Given any compact set
$K\subset \tilde{M}$, $\pi(K)$ is compact in $M$, so by the loops to
infinity property on $M$, there is a curve $\bar{C}\in M\setminus
\pi(K)$
which is homotopic along $\gamma$ to $\pi \circ C$.  We can lift
the homotopy map to the cover $\tilde{M}$, and it is easy to see
that this gives us a curve, $\tilde{\bar{C}}\in \tilde{M}\setminus K$
which is homotopic along $\tilde{\gamma}$ to $C$.

Thus the orientable double cover, $\tilde{M}$, satisfies all the 
conditions
of Proposition~\ref{PartI} and $H_{n-1}(\tilde{M}, G)$ is trivial. 
We claim that this implies that $H_{n-1}(M,Z)$ is trivial as well.

Given any simplicial map $\sigma: \Delta_k \to M$,
there is a continuous lift of $\sigma$,
$\tilde{\sigma}: \Delta_k \to \tilde{M}$,
which is unique up to deck transform.  Thus, there is a
unique chain in $C_k(\tilde{M},\ZZ)$,
\begin{equation}
f(\sigma)= \tilde{\sigma} + g\tilde{\sigma},
\end{equation}
where $g$ is the deck transform acting on the double cover.
We can extend $f$ to a well defined map,
$f: C_k(M, \ZZ) \to C_k(\tilde{M},\ZZ)$.  In
fact, $f$ commutes with the boundary operator, $\partial$,
so $f_*:H_k(M, \ZZ) \to H_k(\tilde{M},\ZZ)$ is well defined.
Note also, that $\pi_* \circ f_*([\sigma])=  2[\sigma]$.

Thus, given any $h \in H_{n-1}(M,\ZZ)$, $f_*(h)$ is in the trivial
group, $H_{n-1}(\tilde{M},\ZZ)$, so
\begin{equation}
2h= \pi_* ( f_*(h))=  \pi_*(0)= 0.
\end{equation}
Thus, $2H_{n-1}(M,\ZZ)=0$ and
\begin{equation}  \label{otimes}
H_{n-1}(M, \ZZ) \otimes \ZZ_2 =  \frac{ 
H_{n-1}(M,\ZZ)}{2H_{n-1}(M,\ZZ)}
=  H_{n-1}(M,\ZZ).
\end{equation}
However, by the Universal Coefficient Theorem (c.f. [Mun, Thm 55.1]),
we have a short exact sequence,
\begin{equation} \label{note1}
0 \to H_{n-1}(M, \ZZ) \otimes \ZZ_2 \to  H_{n-1}(M, \ZZ_2) \to
H_{n-2}(M, \ZZ) \ast \ZZ_2 \to 0.
\end{equation}
Now, by Proposition~\ref{PartI}, $H_{n-1}(M, \ZZ_2)$ is trivial,
so $H_{n-1}(M, \ZZ) \otimes \ZZ_2= 0$.  By (\ref{otimes}),
$H_{n-1}(M, \ZZ)=0$ and the proposition follows.

\qed

Note that in the above proof of Proposition~\ref{PartIII}, we
only use nonnegative Ricci curvature in the first step to eliminate
the possibility of a two-ended double cover.  This step could be
proven without the nonnegative Ricci curvature condition, but it makes
the proof unnecessarily long.  Note also that it is in (\ref{note1}),
when we apply the universal coefficient theorem, that we get our
result for $\ZZ$.  This step does not extend to a result
with the more general $\ZZ_k$ or arbitrary abelian $G$.

\bigskip   \bigskip  
\noindent{\bf Proof of Proposition~\ref{PartIV}}
If $M^n$ does not have a ray satisfying the loops to infinity
property then by
[So2],
$M^n$ is a flat normal bundle over a compact totally
geodesic submanifold, $N^{n-1}$.

Since $M^n$ has only
one end, $N^{n-1}$ is orientable if and only if
$M^n$ is unorientable.  Since the base space of a
normal bundle is homotopic to the bundle,
$H_{n-1}(M, G)=H_{n-1}(N, G)$.  Thus
$$
H_{n-1}(M,G)= H_{n-1}(N, G)=
\cases{G,&if  M is unorientable \cr
 \ker(G \stackrel {\times 2} {\rightarrow} G) ,& if M is orientable }
$$
(c.f. [Mun Cor 65.5]).
\qed

\bigskip

\noindent{\bf Proof of Corollary~\ref{Cor2}:}
If  $M$ is neither a isometrically split manifold nor a flat
normal bundle over a compact totally geodesic submanifold,
then by Theorem~\ref{A}, we know $H_{n-1}(M,G)=0$.
The Universal Coefficient Theorem (c.f. [Mun, Thm 55.1]),
states that there is a short exact sequence,
\begin{equation}
0 \to H_{n-1}(M, \ZZ) \otimes G \to  H_{n-1}(M, G) \to
H_{n-2}(M, \ZZ) \ast G \to 0.
\end{equation}
Substituting,  $H_{n-1}(M,G)=0$, we get
\begin{equation}
H_{n-2}(M, \ZZ) \ast G=0.
\end{equation}
In particular, substituting $G=Z_k$, we have
\begin{equation}
0=H_{n-2}(M, \ZZ) \ast \ZZ_k
= \frac{H_{n-2}(M, \ZZ)}{k H_{n-2}(M, \ZZ)}.
\end{equation}
Thus for any finite number $k$,  $H_{n-2}(M, \ZZ)$ no elements of order 
$k$.  \ProofEnd

\end{document}